\documentstyle[amssymb,12pt]{article}
\newtheorem{th}{Theorem}[section]

\newtheorem{defn}[th]{Definition}
\newenvironment{defn-new}{\begin{defn} \em}{\end{defn}}
\newtheorem{rem}[th]{Remark}
\newenvironment{rem-new}{\begin{rem} \em}{\end{rem}}
\newtheorem{ex}[th]{Example}
\newenvironment{ex-new}{\begin{ex} \em}{\end{ex}}

\makeatletter \@addtoreset{equation}{section} \makeatother

\begin{document}



\begin{center}
{\Large {\bf A new connection in a Riemannian manifold}} \bigskip

Mukut Mani Tripathi\bigskip
\end{center}

\noindent {\bf Abstract.} In a Riemannian manifold, the existence of a new
connection is proved. In particular cases, this connection reduces to
several symmetric, semi-symmetric and quarter-symmetric connections; even
some of them are not introduced so far. We also find formula for curvature
tensor of this new connection. \medskip

 \medskip


\section{Introduction}

Let $\tilde{\nabla}$ be a linear connection in an $n$-dimensional
differentiable manifold $M$. The torsion tensor $\tilde{T}$ and the
curvature tensor $\tilde{R}$ of $\tilde{\nabla}$ are given respectively by
\[
\tilde{T}\left( X,Y\right) \equiv \tilde{\nabla}_{X}Y\ -\ \tilde{\nabla}%
_{Y}X\ -\ [X,Y],
\]
\[
\tilde{R}\left( X,Y\right) Z\equiv \tilde{\nabla}_{X}\tilde{\nabla}_{Y}Z\ -\
\tilde{\nabla}_{Y}\tilde{\nabla}_{X}Z\ -\ \tilde{\nabla}_{[X,Y]}Z.
\]
The connection $\tilde{\nabla}$ is symmetric if its torsion tensor $\tilde{T}
$ vanishes, otherwise it is non-symmetric. The connection $\tilde{\nabla}$
is a metric connection if there is a Riemannian metric $g$ in $M$ such that $%
\tilde{\nabla}g=0$, otherwise it is non-metric. It is well known that a
linear connection is symmetric and metric if and only if it is the
Levi-Civita connection. In 1973, B. G. Schmidt \cite{Schmidt-73} proved that
if the holonomy group of $\tilde{\nabla}$ is a subgroup of the orthogonal
group ${\cal O}(n)$, then $\tilde{\nabla} $ is the Levi-Civita connection of
a Riemannian metric.

In 1930, H. A. Hayden \cite{Hayden-32} introduced a metric connection $%
\tilde{\nabla}$ with a non-zero torsion on a Riemannian manifold. Such a
connection is called a Hayden connection. On the other hand, in a Riemannian
manifold given a $1$-form $\omega $, the {\em Weyl connection} $\tilde{\nabla%
}$ constructed with $\omega $ and its associated vector $B$ (G. B. Folland
1970, \cite{Folland-70}) is a symmetric non-metric connection. In fact, the
Riemannian metric of the manifold is recurrent with respect to the Weyl
connection with the recurrence $1$-form $\omega $, that is, $\tilde{\nabla}%
g=\omega \otimes g$. Another symmetric non-metric connection is projectively
related to the Levi-Civita connection (cf. K. Yano \cite{Yano-70-integral},
D. Smaranda \cite{Smaranda-83}). \medskip

In 1924, A. Friedmann and J. A. Schouten (\cite{Friedmann-Schouten-24}, \cite%
{Schouten-54}) introduced the idea of a semi-symmetric linear connection in
a differentiable manifold. A linear connection is said to be a
semi-symmetric connection if its torsion tensor $\tilde{T}$ is of the form
\begin{equation}
\tilde{T}\left( X,Y\right) =u\left( Y\right) X\ -\ u\left( X\right) Y,
\label{eq-torsion-semi}
\end{equation}
where $u$ is a $1$-form. A Hayden connection with the torsion tensor of the
form (\ref{eq-torsion-semi}) is a semi-symmetric metric connection, which
appeared in a study of E. Pak \cite{Pak-69}. In 1970, K. Yano \cite%
{Yano-70-semi} considered a semi-symmetric metric connection and studied
some of its properties. He proved that a Riemannian manifold is conformally
flat if and only if it admits a semi-symmetric metric connection whose
curvature tensor vanishes identically. He also proved that a Riemannian
manifold is of constant curvature if and only if it admits a semi-symmetric
metric connection for which the manifold is a group manifold, where a group
manifold \cite{Eisenhart-33} is a differentiable manifold admitting a linear
connection $\tilde{\nabla}$ such that its curvature tensor $\tilde{R}$
vanishes and its torsion tensor $\tilde{T}$ is covariantly constant with
respect to $\tilde{\nabla}$. In \cite{Tamassy-Binh-89}, L. Tam\'{a}ssy and
T. Q. Binh proved that if in a Riemannian manifold of dimension $\geq 4$, $%
\tilde{\nabla}$ is a metric linear connection of nonvanishing constant
curvature for which
\[
\tilde{R}\left( X,Y\right) Z+\tilde{R}\left( Y,Z\right) X+\tilde{R}\left(
Z,X\right) Y=0,
\]
then $\tilde{\nabla}$ is the Levi-Civita connection. Some different kind of
semi-symmetric non-metric connections are studied in \cite{Agashe-Chafle-92}%
, \cite{Andonie-Smaranda-77}, \cite{Liang-94} and \cite{Sengupta-De-Binh}.
\medskip

In 1975, S. Golab \cite{Golab-75} defined and studied quarter-symmetric
linear connections in differentiable manifolds. A linear connection is said
to be a quarter-symmetric connection if its torsion tensor $\tilde{T}$ is of
the form
\begin{equation}
\tilde{T}\left( X,Y\right) =u\left( Y\right) \varphi X\ -\ u\left( X\right)
\varphi Y,  \label{eq-torsion-quarter}
\end{equation}
where $u$ is a $1$-form and $\varphi $ is a tensor of type $\left(
1,1\right) $. In 1980, R. S. Mishra and S. N. Pandey \cite{Mishra-Pandey-80}
studied quarter-symmetric metric connections and, in particular, Ricci
quarter-symmetric metric connections.
Note that a quarter-symmetric metric connection is a Hayden connection with
the torsion tensor of the form (\ref{eq-torsion-quarter}). Studies of
various types of quarter-symmetric metric connections and their properties
include \cite{Mishra-Pandey-80}, \cite{Rastogi-78}, \cite{Rastogi-87} and
\cite{Yano-Imai-77} among others.

\medskip Now we ask the following question. Can there be a unified theory of
connections which unifies the concepts of various metric connections such as
a semi-symmetric metric connection and a quarter-symmetric metric connection
and various non-metric connections such as the Weyl connection and different
kind of semi-symmetric non-metric connections?
Surprisingly, we answer to this question in affirmative, and in this paper
we prove the existence of a new connection, which unifies all these
previously known connections and some other connections not introduced so
far. We also find formula for curvature tensor of this new connection. In
the last we list a number of connections as particular cases.

\section{A new connection}

In this section we prove the existence of a new connection in the following

\begin{th}
\label{th-T-connection-existence} Let $M$ be an $n$-dimensional Riemannian
manifold equipped with the Levi-Civita connection $\nabla $ of its
Riemannian metric $g$. Let $f_{1}$, $f_{2}$ be functions in $M$, $u$, $u_{1}$%
, $u_{2}$ are $1$-forms in $M$ and $\varphi $ is a $\left( 1,1\right) $
tensor field in $M$. Let
\begin{equation}
u\left( X\right) \equiv g\left( U,X\right) ,\quad u_{1}\left( X\right)
\equiv g\left( U_{1},X\right) ,\quad u_{2}\left( X\right) \equiv g\left(
U_{2},X\right) ,  \label{eq-defn-U-U1-U2}
\end{equation}
\begin{equation}
g\left( \varphi X,Y\right) \equiv \Phi \left( X,Y\right) =\Phi _{1}\left(
X,Y\right) \ +\ \Phi _{2}\left( X,Y\right) ,  \label{eq-defn-Phi-Phi1-Phi2}
\end{equation}
where $\Phi _{1}$ and $\Phi _{2}$ are symmetric and skew-symmetric parts of
the $\left( 0,2\right) $ tensor $\Phi $ such that
\begin{equation}
\Phi _{1}\left( X,Y\right) \equiv g\left( \varphi _{1}X,Y\right) ,\quad \Phi
_{2}\left( X,Y\right) \equiv g\left( \varphi _{2}X,Y\right) .
\label{eq-rel-Phi-varphi}
\end{equation}
Then there exists a unique connection $\tilde{\nabla}$ in $M$ given by
\begin{eqnarray}
\tilde{\nabla}_{X}Y &=&\nabla _{X}Y\ +\ u\left( Y\right) \varphi _{1}X\ -\
u\left( X\right) \varphi _{2}Y\ -\ g\left( \varphi _{1}X,Y\right) U
\nonumber \\
&&-\ f_{1}\{u_{1}\left( X\right) Y\ +\ u_{1}\left( Y\right) X\ -\ g\left(
X,Y\right) U_{1}\}  \nonumber \\
&&-\ f_{2}\,g\left( X,Y\right) U_{2},  \label{eq-T-connection}
\end{eqnarray}
which satisfies
\begin{equation}
\tilde{T}\left( X,Y\right) =u\left( Y\right) \varphi X\ -\ u\left( X\right)
\varphi Y,  \label{eq-cond-torsion}
\end{equation}
and
\begin{eqnarray}
\left( \tilde{\nabla}_{X}g\right) \left( Y,Z\right) &=&2\,f_{1}u_{1}\left(
X\right) g\left( Y,Z\right)  \nonumber \\
&&\ +\ f_{2}\left\{ \frac {} {}u_{2}\left( Y\right) g\left( X,Z\right) \ +\
u_{2}\left( Z\right) g\left( X,Y\right) \right\} ,  \label{eq-cond-metric}
\end{eqnarray}
where $\tilde{T}$ is the torsion tensor of $\tilde{\nabla}$.
\end{th}

\noindent {\bf Proof.} Let $\tilde{\nabla}$ be a linear connection in $M$
given by
\begin{equation}
\tilde{\nabla}_{X}Y=\nabla _{X}Y\ +\ H\left( X,Y\right) .  \label{eq-cond-1}
\end{equation}
Now, we shall determine the tensor field $H$ such that $\tilde{\nabla}$
satisfies (\ref{eq-cond-torsion}) and (\ref{eq-cond-metric}). From (\ref%
{eq-cond-1}) we have
\begin{equation}
\tilde{T}\left( X,Y\right) =H\left( X,Y\right) \ -\ H\left( Y,X\right) .
\label{eq-cond-2}
\end{equation}
Denote
\begin{equation}
G\left( X,Y,Z\right)\equiv \left( \tilde{\nabla}_{X}g\right) \left(
Y,Z\right) .  \label{eq-G}
\end{equation}
From (\ref{eq-cond-1}) and (\ref{eq-G}), we have 
\begin{equation}
g\left( H\left( X,Y\right) ,Z\right) \ +\ g\left( H\left( X,Z\right)
,Y\right) =\ -\ \,G\left( X,Y,Z\right) .  \label{eq-cond-3}
\end{equation}
From (\ref{eq-cond-2}), (\ref{eq-cond-3}), (\ref{eq-G}) and (\ref%
{eq-cond-metric}), we have
\begin{eqnarray*}
&&g\left( \tilde{T}\left( X,Y\right) ,Z\right) \ +\ g\left( \tilde{T}\left(
Z,X\right) ,Y\right) \ +\ g\left( \tilde{T}\left( Z,Y\right) ,X\right) \\
&=&g\left( H\left( X,Y\right) ,Z\right) \ -\ g\left( H\left( Y,X\right)
,Z\right) \ +\ g\left( H\left( Z,X\right) ,Y\right) \\
&&\ -\ g\left( H\left( X,Z\right) ,Y\right) \ +\ g\left( H\left( Z,Y\right)
,X\right) \ -\ g\left( H\left( Y,Z\right) ,X\right) \\
&=&g\left( H\left( X,Y\right) ,Z\right) \ -\ g\left( H\left( X,Z\right)
,Y\right) \ +\ G\left( Y,X,Z\right) \ -\ G\left( Z,X,Y\right) \\
&=&2g\left( H\left( X,Y\right) ,Z\right) \ +\ G\left( X,Y,Z\right) \ +\
G\left( Y,X,Z\right) \ -\ G\left( Z,X,Y\right) \\
&=&2g\left( H\left( X,Y\right) ,Z\right) \\
&&\ +\ \,\left\{ 2f_{1}u_{1}\left( X\right) g\left( Y,Z\right) \ +\
f_{2}u_{2}\left( Y\right) g\left( X,Z\right) \ +\ f_{2}u_{2}\left( Z\right)
g\left( X,Y\right) \right\} \\
&&\ +\ \,\left\{ 2f_{1}u_{1}\left( Y\right) g\left( X,Z\right) \ +\
f_{2}u_{2}\left( X\right) g\left( Y,Z\right) \ +\ f_{2}u_{2}\left( Z\right)
g\left( X,Y\right) \right\} \\
&&\ -\ \,\left\{ 2f_{1}u_{1}\left( Z\right) g\left( X,Y\right) \ +\
f_{2}u_{2}\left( X\right) g\left( Y,Z\right) \ +\ f_{2}u_{2}\left( Y\right)
g\left( X,Z\right) \right\} ,
\end{eqnarray*}
which in view of (\ref{eq-defn-U-U1-U2}) implies that
\begin{eqnarray}
H\left( X,Y\right) &=&\frac{1}{2}\{ \tilde{T}\left( X,Y\right) \ +\ \tilde{T}%
^{\prime }\left( X,Y\right) \ +\ \tilde{T}^{\prime }\left( Y,X\right) \}
\nonumber \\
&&\ -\ f_{1}\left\{ u_{1}\left( X\right) Y\ +\ u_{1}\left( Y\right) X\ -\
f_{1}g\left( X,Y\right) U_{1}\right\}  \nonumber \\
&&\ -\ f_{2}\,g\left( X,Y\right) U_{2}\!,  \label{eq-cond-4}
\end{eqnarray}
where
\begin{equation}
g(\tilde{T}^{\prime }\left( X,Y\right) ,Z)=g\left( \tilde{T}\left(
Z,X\right) ,Y\right) .  \label{eq-cond-5}
\end{equation}
Using (\ref{eq-cond-torsion}), (\ref{eq-defn-U-U1-U2}), (\ref%
{eq-defn-Phi-Phi1-Phi2}) and (\ref{eq-rel-Phi-varphi}) in (\ref{eq-cond-5}),
we get
\begin{eqnarray*}
g(\tilde{T}^{\prime }\left( X,Y\right) ,Z) &=&g\left( u\left( X\right)
\varphi Z\ -\ u\left( Z\right) \varphi X,Y\right) \\
&=&u\left( X\right) \Phi \left( Z,Y\right) \ -\ g\left( U,Z\right) \Phi
\left( X,Y\right) \\
&=&u\left( X\right) g\left( \varphi _{1}Z,Y\right) \ +\ u\left( X\right)
g\left( \varphi _{2}Z,Y\right) \\
&&\ -\ g\left( U,Z\right) \Phi \left( X,Y\right) ,
\end{eqnarray*}
which implies that
\begin{equation}
\tilde{T}^{\prime }\left( X,Y\right) =u\left( X\right) \varphi _{1}Y\ -\
u\left( X\right) \varphi _{2}Y\ -\ \Phi \left( X,Y\right) U.
\label{eq-cond-6}
\end{equation}
In view of (\ref{eq-cond-torsion}), (\ref{eq-cond-4}) and (\ref{eq-cond-6}),
we get
\begin{eqnarray}
H\left( X,Y\right) &=&u\left( Y\right) \varphi _{1}X\ -\ u\left( X\right)
\varphi _{2}Y\ -\ g\left( \varphi _{1}X,Y\right) U  \nonumber \\
&&\ -\ f_{1}\left\{ u_{1}\left( X\right) Y\ +\ u_{1}\left( Y\right) X\ -\
g\left( X,Y\right) U_{1}\right\}  \nonumber \\
&&\ -\ f_{2}\,g\left( X,Y\right) U_{2}  \label{eq-H(X,Y)}
\end{eqnarray}
and hence, $\tilde{\nabla}$ is given by (\ref{eq-T-connection}). $\medskip $

Conversely, a connection given by (\ref{eq-T-connection}) satisfies the
conditions (\ref{eq-cond-torsion}) and (\ref{eq-cond-metric}). $\blacksquare
$

\section{Curvature tensor}

The curvature tensor $\tilde{R}$ of the connection $\tilde{\nabla}$ is given
by
\[
\tilde{R}\left( X,Y\right) Z=\tilde{\nabla}_{X}\tilde{\nabla}_{Y}Z\ -\
\tilde{\nabla}_{Y}\tilde{\nabla}_{X}Z\ -\ \tilde{\nabla}_{\left[ X,Y\right]
}Z.
\]
For a function $h$, in view of (\ref{eq-T-connection}) we have
\begin{eqnarray}
\tilde{\nabla}_{X}\left( hY\right) &=&\left( Xh\right) Y\ +\ h\nabla _{X}Y
\nonumber \\
&&\ +\ hu\left( Y\right) \varphi _{1}X\ -\ hu\left( X\right) \varphi _{2}Y\
-\ hg\left( \varphi _{1}X,Y\right) U  \nonumber \\
&&\ -\ f_{1}\left\{ hu_{1}\left( X\right) Y\ +\ hu_{1}\left( Y\right) X\ -\
hg\left( X,Y\right) U_{1}\right\}  \nonumber \\
&&\ -\ f_{2}h\,g\left( X,Y\right) U_{2}.  \label{eq-delbar-X-(h-psi-Y)}
\end{eqnarray}
Also from (\ref{eq-T-connection}) it follows that
\begin{eqnarray}
\tilde{\nabla}_{X}\tilde{\nabla}_{Y}Z\ &=&\tilde{\nabla}_{X}\nabla _{Y}Z\ +\
\tilde{\nabla}_{X}\left( u\left( Z\right) \varphi _{1}Y\right) \ -\ \tilde{%
\nabla}_{X}\left( u\left( Y\right) \varphi _{2}Z\right)  \nonumber \\
&&\ -\ \tilde{\nabla}_{X}\left( g\left( \varphi _{1}Y,Z\right) U\right) \ -\
\tilde{\nabla}_{X}\left( f_{1}u_{1}\left( Y\right) Z\right)  \nonumber \\
&&\ -\ \tilde{\nabla}_{X}\left( f_{1}u_{1}\left( Z\right) Y\right) \ +\
\tilde{\nabla}_{X}\left( f_{1}g\left( Y,Z\right) U_{1}\right)  \nonumber \\
&&\ -\ \tilde{\nabla}_{X}\left( f_{2}\,g\left( Y,Z\right) U_{2}\right) .
\label{eq-delbar-X-delbar-Y-Z}
\end{eqnarray}
We shall need the following definitions. Let $\eta $ be a $1$-form and $\xi$
be its associated vector field such that
\[
\eta \left( X\right) =g\left( \xi ,X\right) .
\]
We define
\begin{equation}
\beta(\eta ,X,Y) = (\nabla_{X}\eta)Y + u(X)\eta(\varphi_{2}Y) - \eta(\varphi
_{1}X)u(Y) + \eta(U)g(\varphi_{1}X,Y),  \label{eq-beta(eta,X,Y)}
\end{equation}
\begin{equation}
B\left( \eta ,X\right) =\nabla _{X}\xi -u\left( X\right) \varphi _{2}\xi
-\eta \left( \varphi _{1}X\right) U+\eta \left( U\right) \varphi _{1}X,
\label{eq-B(eta,X)}
\end{equation}
such that
\begin{equation}
\beta \left( \eta ,X,Y\right) =g\left( B\left( \eta ,X\right) ,Y\right) .
\label{eq-alpha-A}
\end{equation}
Let
\begin{equation}
\alpha \left( \eta ,X,Y\right) =\beta \left( \eta ,X,Y\right) -\frac{1}{2}\,
\eta \left( U\right) g\left( \varphi _{1}X,Y\right) .  \label{eq-beta-alpha}
\end{equation}
Then
\begin{equation}
\alpha \left( \eta ,X,Y\right) =g\left( A\left( \eta ,X\right) ,Y\right) ,
\label{eq-beta-B}
\end{equation}
where
\begin{equation}
A\left( \eta ,X\right) =B\left( \eta ,X\right) -\frac{1}{2}\,\eta \left(
U\right) \varphi _{1}X.  \label{eq-B-A}
\end{equation}
Explicitly, we have
\begin{equation}
\alpha(\eta,X,Y) = (\nabla_{X}\eta)Y + u(X)\eta(\varphi_{2}Y) -
\eta(\varphi_{1}X)u(Y) + \frac{1}{2}\,\eta(U)g(\varphi _{1}X,Y),
\label{eq-alpha(eta,X,Y)}
\end{equation}
\begin{equation}
A\left( \eta ,X\right) =\nabla _{X}\xi -u\left( X\right) \varphi _{2}\xi
-\eta \left( \varphi _{1}X\right) U+\frac{1}{2}\eta \left( U\right) \varphi
_{1}X,  \label{eq-A(eta,X)}
\end{equation}
We also define
\begin{equation}
\mu \left( X,Y\right) =\left( \nabla _{X}\varphi _{1}\right) Y\ -\ u\left(
X\right) \varphi _{2}\varphi _{1}Y,  \label{eq-mu}
\end{equation}
\begin{equation}
R_{0}\left( X,Y\right) Z=g\left( Y,Z\right) X-g\left( X,Z\right) Y.
\label{eq-R-0}
\end{equation}

In view of (\ref{eq-T-connection}), (\ref{eq-delbar-X-(h-psi-Y)}) and (\ref%
{eq-delbar-X-delbar-Y-Z}) we obtain the following formula for curvature $%
\tilde{R}$ of the connection $\tilde{\nabla}$
\begin{eqnarray}
\tilde{R}\left( X,Y\right) Z &=&R\left( X,Y\right) Z\ -\ 2du\left(
X,Y\right) \varphi _{2}Z  \nonumber \\
&&\left. -\ \alpha \left( u,Y,Z\right) \varphi _{1}X\ +\ \alpha \left(
u,X,Z\right) \varphi _{1}Y\right.  \nonumber \\
&&\left. -\ g\left( \varphi _{1}Y,Z\right) A\left( u,X\right) \ +\ g\left(
\varphi _{1}X,Z\right) A\left( u,Y\right) \right.  \nonumber \\
&&\left. -\ R_{0}\left( U,\mu \left( X,Y\right) -\mu \left( Y,X\right)
\right) Z\right.  \nonumber \\
&&\left. +\ u\left( X\right) \left( \nabla _{Y}\varphi _{2}\right) Z\ -\
u\left( Y\right) \left( \nabla _{X}\varphi _{2}\right) Z\right.  \nonumber \\
&&-\ f_{1}\left\{ 2du_{1}\left( X,Y\right) Z-\beta \left( u_{1},Y,Z\right)
X+\beta \left( u_{1},X,Z\right) Y\right.  \nonumber \\
&&\qquad \quad \left. -\ g\left( Y,Z\right) B\left( u_{1},X\right) \ +\
g\left( X,Z\right) B\left( u_{1},Y\right) \right.  \nonumber \\
&&\qquad \quad \left. +\ u\left( Y\right) R_{0}\left( \varphi X,U_{1}\right)
Z\ -\ u\left( X\right) R_{0}\left( \varphi Y,U_{1}\right) Z\right\}
\nonumber \\
&&+\ f_{2}\left\{ \,g\left( \varphi Y,Z\right) u\left( X\right) U_{2}\ -\
\,g\left( \varphi X,Z\right) u\left( Y\right) U_{2}\right.  \nonumber \\
&&\qquad \quad \left. -\ g\left( Y,Z\right) B\left( u_{2},X\right) \ +\
g\left( X,Z\right) B\left( u_{2},Y\right) \right\}  \nonumber \\
&&-\ \left( f_{1}\right) ^{2}\left\{ g\left( Y,Z\right) R_{0}\left(
X,U_{1}\right) U_{1}-g\left( X,Z\right) R_{0}\left( Y,U_{1}\right)
U_{1}\right.  \nonumber \\
&&\qquad \qquad \left. -\ u_{1}\left( Z\right) R_{0}\left( X,Y\right)
U_{1}\right\}  \nonumber \\
&&+\left( f_{2}\right) ^{2}\left\{ g\left( Y,Z\right) u_{2}\left( X\right)
U_{2}-g\left( X,Z\right) u_{2}\left( Y\right) U_{2}\right\}  \nonumber \\
&&+\ f_{1}f_{2}\left\{ g\left( Y,Z\right) \left( R_{0}\left( X,U_{2}\right)
U_{1}\ -\ u_{2}\left( X\right) U_{1}\right) \right.  \nonumber \\
&&\qquad \qquad \left. \ -\ g\left( X,Z\right) \left( R_{0}\left(
Y,U_{2}\right) U_{1}\ -\ u_{2}\left( Y\right) U_{1}\right) \right\}
\nonumber \\
&&-\ \left( Xf_{1}\right) \left\{ u_{1}\left( Y\right) Z\ +\ u_{1}\left(
Z\right) Y\ -\ g\left( Y,Z\right) U_{1}\right\}  \nonumber \\
&&+\ \left( Yf_{1}\right) \left\{ u_{1}\left( X\right) Z\ +\ u_{1}\left(
Z\right) X\ -\ g\left( X,Z\right) U_{1}\right\}  \nonumber \\
&&-\ \left( Xf_{2}\right) g\left( Y,Z\right) U_{2}\ +\ \left( Yf_{2}\right)
g\left( X,Z\right) U_{2},  \label{eq-curvature-formula}
\end{eqnarray}
where (\ref{eq-alpha(eta,X,Y)}), (\ref{eq-A(eta,X)}), (\ref{eq-beta(eta,X,Y)}%
), (\ref{eq-B(eta,X)}), (\ref{eq-mu}) and (\ref{eq-R-0}) are used.

\section{Particular cases}

In this section, we list the following seventeen particular cases.

\subsection{Quarter-symmetric metric connections}

\begin{enumerate}
\item[{\bf (1)}] $f_{1}=0=f_{2}$. Then we obtain a {\bf quarter-symmetric
metric connection} $\tilde{\nabla}$ given by (K. Yano and T. Imai \cite%
{Yano-Imai-82}, eq. (3.3))
\[
\tilde{\nabla}_{X}Y=\nabla _{X}Y\ +\ u\left( Y\right) \varphi _{1}X\ -\
u\left( X\right) \varphi _{2}Y\ -\ g\left( \varphi _{1}X,Y\right) U.
\]

\item[{\bf (2)}] $f_{1}=0=f_{2}$, $\varphi _{2}=0$. In this case $\varphi
=\varphi _{1}$ and $g\left( \varphi X,Y\right) =g\left( X,\varphi Y\right) $%
. Then we obtain a quarter-symmetric metric connection $\tilde{\nabla}$
given by (R. S. Mishra and S. N. Pandey \cite{Mishra-Pandey-80}, eq. (1.6))
\[
\tilde{\nabla}_{X}Y=\nabla _{X}Y\ +\ u\left( Y\right) \varphi X\ -\
g(\varphi X,Y)U.
\]
In particular, if $f_{1}=0=f_{2}$, $\varphi =Q$, where $Q$ is the Ricci
operator, then we obtain the {\em Ricci quarter-symmetric metric connection}
$\tilde{\nabla}$ given by (R. S. Mishra and S. N. Pandey \cite%
{Mishra-Pandey-80}, eq. (2.2))
\[
\tilde{\nabla}_{X}Y=\nabla _{X}Y\ +\ u\left( Y\right) QX\ -\ S\left(
X,Y\right) U,
\]
where $S$ is the Ricci tensor. The torsion of this connection satisfies
\[
\tilde{T}\left( X,Y\right) =u\left( Y\right) QX\ -\ u\left( X\right) QY.
\]

\item[{\bf (3)}] $f_{1}=0=f_{2}$, $\varphi _{1}=0$. Then we obtain a
quarter-symmetric metric connection $\tilde{\nabla}$ given by (K. Yano and
T. Imai \cite{Yano-Imai-82}, eq. (3.6))
\[
\tilde{\nabla}_{X}Y=\nabla _{X}Y\ -\ u\left( X\right) \varphi Y.
\]
This connection was also introduced by R. S. Mishra and S. N. Pandey in an
almost Hermitian manifold. (cf. \cite{Mishra-Pandey-80}, eq. (3.1)). In \cite%
{Ojha-Prasad-74}, R. H. Ojha and S. Prasad defined this type of connection
in an almost contact metric manifold and called it a {\em semi-symmetric
metric} $S${\em -connection}.
\end{enumerate}

\subsection{Quarter-symmetric non-metric connections}

\begin{enumerate}
\item[{\bf (4)}] $f_{1}\neq 0$, $f_{2}=0$, $\varphi _{2}=0$. Then we obtain
a {\em quarter-symmetric recurrent-metric connection} $\tilde{\nabla}$ given
by
\begin{eqnarray*}
\tilde{\nabla}_{X}Y &=&\nabla _{X}Y\ +\ u\left( Y\right) \varphi X\ -\
g\left( \varphi X,Y\right) U \\
&&\ -\ f_{1}\left\{ u_{1}\left( X\right) Y\ +\ u_{1}\left( Y\right) X\ -\
g\left( X,Y\right) U_{1}\right\} .
\end{eqnarray*}
This connection satisfies $\tilde{\nabla}g=2f_{1}u_{1}\otimes g$.

\item[{\bf (5)}] $f_{1}=1$, $f_{2}=0$, $\varphi _{2}=0$, $u_{1}=u$. Then we
obtain a {\em special quarter-symmetric recurrent-metric connection} $\tilde{%
\nabla}$ given by
\begin{eqnarray*}
\tilde{\nabla}_{X}Y &=&\nabla _{X}Y\ +\ u\left( Y\right) \varphi X\ -\
g\left( \varphi X,Y\right) U \\
&&\ -\ u\left( X\right) Y\ -\ u\left( Y\right) X\ +\ g\left( X,Y\right) U.
\end{eqnarray*}
This connection satisfies $\tilde{\nabla}g=2u\otimes g$.

\item[{\bf (6)}] $f_{1}\neq 0$, $f_{2}=0$, $\varphi _{1}=0$. Then we obtain
a {\em quarter-symmetric recurrent-metric connection} $\tilde{\nabla}$ given
by
\begin{eqnarray*}
\tilde{\nabla}_{X}Y &=&\nabla _{X}Y\ -\ u\left( X\right) \varphi Y \\
&&\ -\ f_{1}\left\{ u_{1}\left( X\right) Y\ +\ u_{1}\left( Y\right) X\ -\
g\left( X,Y\right) U_{1}\right\} .
\end{eqnarray*}
This connection satisfies $\tilde{\nabla}g=f_{1}u_{1}\otimes g$.

\item[{\bf (7)}] $f_{1}=1$, $f_{2}=0$, $\varphi _{1}=0$, $u_{1}=u$. Then we
obtain a {\em special quarter-symmetric recurrent-metric connection} $\tilde{%
\nabla}$ given by
\[
\tilde{\nabla}_{X}Y=\nabla _{X}Y\ -\ u\left( X\right) \varphi Y\ -\ u\left(
X\right) Y\ -\ u\left( Y\right) X\ +\ g\left( X,Y\right) U.
\]
This connection satisfies $\tilde{\nabla}g=2u\otimes g$.

\item[{\bf (8)}] $f_{1}=0$, $f_{2}\neq 0$, $\varphi _{2}=0$. Then we obtain
a {\em quarter-symmetric non-metric connection} $\tilde{\nabla}$ given by
\[
\tilde{\nabla}_{X}Y=\nabla _{X}Y\ +\ u\left( Y\right) \varphi X\ -\ g\left(
\varphi X,Y\right) U\ -\ f_{2}\,g\left( X,Y\right) U_{2}.
\]
This connection satisfies
\[
\left( \tilde{\nabla}_{X}g\right) \left( Y,Z\right) =f_{2}\left\{
u_{2}\left( Y\right) g\left( X,Z\right) \ +\ u_{2}\left( Z\right) g\left(
X,Y\right) \right\} .
\]

\item[{\bf (9)}] $f_{1}=0$, $f_{2}\neq 0$, $\varphi _{2}=0$, $u_{2}=u$. Then
we obtain a {\em quarter-symmetric non-metric connection} $\tilde{\nabla}$
given by
\[
\tilde{\nabla}_{X}Y=\nabla _{X}Y\ +\ u\left( Y\right) \varphi X\ -\ g\left(
\varphi X,Y\right) U\ -\ f_{2}\,g\left( X,Y\right) U.
\]
This connection satisfies
\[
\left( \tilde{\nabla}_{X}g\right) \left( Y,Z\right) =f_{2}\left\{ u\left(
Y\right) g\left( X,Z\right) \ +\ u\left( Z\right) g\left( X,Y\right)
\right\} .
\]

\item[{\bf (10)}] $f_{1}=0$, $f_{2}\neq 0$, $\varphi _{1}=0$. Then we obtain
a {\em quarter-symmetric non-metric connection} $\tilde{\nabla}$ given by
\[
\tilde{\nabla}_{X}Y=\nabla _{X}Y\ -\ u\left( X\right) \varphi Y\ -\
f_{2}\,g\left( X,Y\right) U_{2}.
\]
This connection satisfies
\[
\left( \tilde{\nabla}_{X}g\right) \left( Y,Z\right) =f_{2}\left\{
u_{2}\left( Y\right) g\left( X,Z\right) \ +\ u_{2}\left( Z\right) g\left(
X,Y\right) \right\} .
\]

\item[{\bf (11)}] $f_{1}=0$, $f_{2}\neq 0$, $\varphi _{1}=0$, $u_{2}=u$.
Then we obtain a {\em quarter-symmetric non-metric connection} $\tilde{\nabla%
}$ given by
\[
\tilde{\nabla}_{X}Y=\nabla _{X}Y\ -\ u\left( X\right) \varphi Y\ -\
f_{2}\,g\left( X,Y\right) U.
\]
This connection satisfies
\[
\left( \tilde{\nabla}_{X}g\right) \left( Y,Z\right) =f_{2}\left\{ u\left(
Y\right) g\left( X,Z\right) \ +\ u\left( Z\right) g\left( X,Y\right)
\right\} .
\]
\end{enumerate}

\subsection{Semi-symmetric metric connection{\bf \ }}

\begin{enumerate}
\item[{\bf (12)}] $f_{1}=0=f_{2}$, $\varphi ={Id}$. Then we obtain a {\em %
semi-symmetric metric connection} $\tilde{\nabla}$ given by (K. Yano, 1970
\cite{Yano-70-semi})
\[
\tilde{\nabla}_{X}Y=\nabla _{X}Y\ +\ u\left( Y\right) X\ -\ g\left(
X,Y\right) U.
\]
\end{enumerate}

\subsection{Semi-symmetric non-metric connections}

\begin{enumerate}
\item[{\bf (13)}] $f_{1}\neq 0$, $f_{2}=0$, $\varphi ={Id}$. Then we obtain
a {\em semi-symmetric recurrent-metric connection} $\tilde{\nabla}$ given by
\begin{eqnarray*}
\tilde{\nabla}_{X}Y &=&\nabla _{X}Y\ +\ u\left( Y\right) X\ -\ g\left(
X,Y\right) U \\
&&\ -\ f_{1}\left\{ u_{1}\left( X\right) Y\ +\ u_{1}\left( Y\right) X\ -\
g\left( X,Y\right) U_{1}\right\} .
\end{eqnarray*}
This connection satisfies $\tilde{\nabla}g=f_{1}u_{1}\otimes g$. In
particular, if $f_{1}=1$, $f_{2}=0$, $\varphi ={Id}$, then we obtain a {\em %
semi-symmetric recurrent-metric connection} $\tilde{\nabla}$ given by (O. C.
Andonie, D. Smaranda 1977, \cite{Andonie-Smaranda-77}; Y. Liang 1994, \cite%
{Liang-94})
\begin{eqnarray*}
\tilde{\nabla}_{X}Y &=&\nabla _{X}Y\ +\ u\left( Y\right) X\ -\ g\left(
X,Y\right) U \\
&&\ -\ u_{1}\left( X\right) Y\ -\ u_{1}\left( Y\right) X\ +\ g\left(
X,Y\right) U_{1},
\end{eqnarray*}
This connection satisfies $\tilde{\nabla}g=2u_{1}\otimes g$. In particular,
if $f_{1}=1$, $f_{2}=0$, $\varphi ={Id}$, $u_{1}=u$, then we obtain a {\em %
semi-symmetric recurrent-metric connection} $\tilde{\nabla}$ given by
\[
\tilde{\nabla}_{X}Y=\nabla _{X}Y\ -\ u\left( X\right) Y.
\]
This connection satisfies $\tilde{\nabla}g=2u\otimes g$.

\item[{\bf (14)}] $f_{1}=0$, $f_{2}\neq 0$, $\varphi ={Id}$. Then we obtain
a {\em semi-symmetric non-metric connection} $\tilde{\nabla}$ given by
\[
\tilde{\nabla}_{X}Y=\nabla _{X}Y\ +\ u\left( Y\right) X\ -\ g\left(
X,Y\right) U\ -\ f_{2}\,g\left( X,Y\right) U_{2}.
\]
This connection satisfies
\[
\left( \tilde{\nabla}_{X}g\right) \left( Y,Z\right) =f_{2}\left\{
u_{2}\left( Y\right) g\left( X,Z\right) \ +\ u_{2}\left( Z\right) g\left(
X,Y\right) \right\} .
\]
In particular, if $f_{1}=0$, $f_{2}=\ -\ 1$, $\varphi ={Id}$, then we obtain
a {\em semi-symmetric non-metric connection} $\tilde{\nabla}$ given by (J.
Sengupta, U. C. De, T. Q. Binh 2000, \cite{Sengupta-De-Binh})
\[
\tilde{\nabla}_{X}Y=\nabla _{X}Y\ +\ u\left( Y\right) X\ -\ g\left(
X,Y\right) U\ +\ \,g\left( X,Y\right) U_{2}.
\]
This connection satisfies
\[
\left( \tilde{\nabla}_{X}g\right) \left( Y,Z\right) =\ -\ u_{2}\left(
Y\right) g\left( X,Z\right) \ -\ u_{2}\left( Z\right) g\left( X,Y\right) .
\]
In particular, if $f_{1}=0$, $f_{2}=\ -\ 1$, $\varphi ={Id}$, $u_{2}=u$,
then we obtain a {\em semi-symmetric non-metric connection} $\tilde{\nabla}$
given by (N. S. Agashe, M. R. Chafle 1992, \cite{Agashe-Chafle-92})
\[
\tilde{\nabla}_{X}Y=\nabla _{X}Y\ +\ u\left( Y\right) X.
\]
This connection satisfies
\[
\left( \tilde{\nabla}_{X}g\right) \left( Y,Z\right) =\ -\ u\left( Y\right)
g\left( X,Z\right) \ -\ u\left( Z\right) g\left( X,Y\right) .
\]
\end{enumerate}

\subsection{Symmetric connections}

\begin{enumerate}
\item[{\bf (15)}] $u=0$. Then we obtain a {\em symmetric non-metric
connection} $\tilde{\nabla}$ given by
\begin{eqnarray*}
\tilde{\nabla}_{X}Y &=&\nabla _{X}Y\ -\ f_{1}\left\{ u_{1}\left( X\right) Y\
+\ u_{1}\left( Y\right) X\ -\ g\left( X,Y\right) U_{1}\right\} \\
&&\ -\ \,f_{2}\,g\left( X,Y\right) U_{2}.
\end{eqnarray*}
This connection satisfies
\begin{eqnarray*}
\left( \tilde{\nabla}_{X}g\right) \left( Y,Z\right) &=&2f_{1}u_{1}\left(
X\right) g\left( Y,Z\right) \\
&&\ +\ f_{2}\left\{ u_{2}\left( Y\right) g\left( X,Z\right) \ +\ u_{2}\left(
Z\right) g\left( X,Y\right) \right\} .
\end{eqnarray*}

\item[{\bf (16)}] $u=0$, $f_{1}=1/2$, $f_{2}=0$, $u_{1}=\omega $, $U_{1}=B$.
Then we obtain a {\em Weyl connection }constructed with $\omega $ and $B$
(G. B. Folland 1970, \cite{Folland-70}) and given by
\[
\tilde{\nabla}_{X}Y=\nabla _{X}Y\ -\ \frac{1}{2}\left\{ \omega \left(
X\right) Y\ +\ \omega \left( Y\right) X\ -\ g\left( X,Y\right) B\right\} .
\]
This connection is a symmetric recurrent-metric connection. The Riemannian
metric $g$ is recurrent with respect to the connection $\tilde{\nabla}$ with
the recurrence $1$-form $\omega $, that is, $\tilde{\nabla}g=\omega \otimes
g $.

\item[{\bf (17)}] $u=0$, $f_{1}=\ -\ 1$, $f_{2}=\ -\ 1$, $u_{1}=u_{2}=\omega
$. In this case, we obtain a {\em symmetric non-metric connection} $\tilde{%
\nabla}$ given by (K. Yano, \cite{Yano-70-integral} eq. (4.5) p. 17; (see
also D. Smaranda 1983, \cite{Smaranda-83}))
\[
\tilde{\nabla}_{X}Y=\nabla _{X}Y\ +\ \omega \left( X\right) Y\ +\ \omega
\left( Y\right) X.
\]
This connection satisfies
\[
\left( \tilde{\nabla}_{X}g\right) \left( Y,Z\right) =\ -\ 2\omega \left(
X\right) g\left( Y,Z\right) \ -\ \omega \left( Y\right) g\left( X,Z\right) \
-\ \omega \left( Z\right) g\left( X,Y\right) .
\]
This connection is projectively related to the Levi-Civita connection $%
\nabla $. Then from (\ref{eq-curvature-formula}) we get 
\begin{eqnarray*}
\tilde{R}\left( X,Y\right) Z &=&R\left( X,Y\right) Z+s\left( X,Z\right)
Y-s\left( Y,Z\right) X \\
&&+\left\{ s\left( X,Y\right) -s\left( Y,X\right) \right\} Z,
\end{eqnarray*}
where
\[
s\left( X,Y\right) =\left( \nabla _{X}\omega \right) Y-\omega \left(
X\right) \omega \left( Y\right) =\left( \tilde{\nabla}_{X}\omega \right) Y.
\]
Note that
\[
s\left( X,Y\right) -s\left( Y,X\right) =2d\omega \left( X,Y\right) ,
\]
that is, $s$ is symmetric if and only if the $1$-form $\omega $ is closed.
\end{enumerate}


\noindent Department of Mathematics and Astronomy

\noindent Lucknow University

\noindent Lucknow-226 007, India

\end{document}